\documentclass[12pt]{amsart}

\usepackage{ulem}
\usepackage[margin=1in]{geometry}
\usepackage{amscd}
\usepackage{amssymb}
\usepackage{amsmath}
\usepackage{amsthm}
\usepackage{enumerate}
\usepackage{tikz}
\usepackage{wrapfig, framed, caption}
\usepackage{float}
\usetikzlibrary{arrows}
\usepackage[font=small,labelfont=bf]{caption}
\usepackage{xcolor}
\usepackage{mathtools}
\usepackage[colorlinks=true, linkcolor=blue, citecolor=blue,
pagebackref=true]{hyperref}
\usepackage{ulem}
\allowdisplaybreaks
\usepackage{mathabx}

{\theoremstyle{definition}
\newtheorem{definition}{Definition}[section]} 
\newtheorem{theorem}[definition]{Theorem}
\newenvironment{theorem*}[1]{{\bf Theorem #1} \begin{itshape}}{\end{itshape}}

\newtheorem{corollary}[definition]{Corollary}
\newenvironment{corollary*}[1]{{\bf Corollary #1} \begin{itshape}}{\end{itshape}}
\newtheorem{proposition}[definition]{Proposition}
\newenvironment{proposition*}[1]{{\bf Proposition #1} \begin{itshape}}{\end{itshape}}

\theoremstyle{remark}
\newtheorem{remark}[definition]{Remark}

\numberwithin{equation}{section}

\renewcommand{\Bbb}[1]{\mathbb{#1}}

\newcommand{\N}{{\Bbb N}}         

\newcommand{\Q}{{\Bbb Q}}         
\newcommand{\R}{{\Bbb R}}        

\newcommand{\Z}{{\Bbb Z}}         

\newcommand{\cB}{{\mathcal B}}
\newcommand{\cC}{{\mathcal C}}

\newcommand{\cH}{{\mathcal H}}

\newcommand{\cL}{{\mathcal L}}
\newcommand{\cM}{{\mathcal M}}

\newcommand{\cR}{{\mathcal R}}
\newcommand{\cS}{{\mathcal S}}

\newcommand{\cU}{{\mathcal U}}

\newcommand{\bx}{\mathbf{x}}
\newcommand{\by}{\mathbf{y}}
\newcommand{\bz}{\mathbf{z}}
\newcommand{\bp}{\mathbf{p}}
\newcommand{\bq}{\mathbf{q}}

\newcommand{\dist}{\operatorname{dist}}
\newcommand{\codim}{\operatorname{codim}}

\newcommand{\bad}{\mathfrak{Bad}}
\newcommand{\Bad}{\mathbf{Bad}}

\DeclareMathOperator{\por}{por}

\DeclareMathOperator{\dimh}{\dim_H}
\DeclareMathOperator{\dimp}{\dim_P}
\DeclareMathOperator{\codimh}{\codim_H}
\DeclareMathOperator{\codimp}{\codim_P}

\title{Porosity in Diophantine approximation}

\author{Henna Koivusalo}
\author{Benjamin Ward}
\date{}
\begin{document}
\frenchspacing
\maketitle
\begin{abstract}
We provide a general framework for obtaining upper bounds on the packing dimension of certain sets arising from Diophantine approximation using the notion of porosity from geometric measure theory. In particular, the framework relies on showing that such sets are mean porous, as introduced by Koskela and Rohde (Math. Ann., 1997). As applications of our framework we show that:
\begin{itemize}
\item For $\varepsilon>0$ sufficiently small, the set of $d$-dimensional $\varepsilon$-badly approximable points has packing codimension comparable to $\varepsilon^{d}$ up to multiplicative constant.
\item Any $C^{2}$ submanifold of $\mathbb{R}^{d}$ with bounded first and second partial derivatives intersected with the set of $d$-dimensional $\varepsilon$-badly approximable points has packing dimension strictly less than the dimension of the submanifold.
\item For any matrix $A\in \mathbb{R}^{m\times n}$ satisfying certain Diophantine conditions, the set of points that avoid the $\varepsilon \|\mathbf{q}\|_{\infty}^{-\frac{n}{m}}$-neighbourhood of the orbit $\{A\mathbf{q} \mod \mathbb{Z}^{m}:\mathbf{q} \in \mathbb{Z}^{n}\}$ has packing codimension greater than $C_{A} \varepsilon^{m}$ for $C_{A}$ a constant depending on $A$ and $m$ only.
\end{itemize}
\end{abstract}
\section{Introduction}
The notion of porosity captures how many holes appear in a set at all scales. The larger the porosity the larger the holes are. See \cite[\S 3.8]{Falconerbook} for a gentle introduction and \cite{ShmerkinSurvey,JarvenpaaSurvey} for surveys. Formally, let $E\subset \R^{d}$, and consider any $x\in E$ and $r>0$. Define
\begin{equation*}
    \por (E,x,r):= \sup\left\{\alpha\geq 0: \text{ there exists $y\in \R^{d}$ such that   } \, \, B(y,\alpha r) \subset B(x,r)\setminus E \right\}.
\end{equation*}
Define the porosity of $E$ at $x$ to be the lower limit of the sequence as $r$ tends to zero. That is,
\begin{equation*}
    \por (E,x) := \liminf_{r\to 0} \por(E,x,r),
\end{equation*}
and the porosity of $E$ to be the infimum over all $x \in E$;
\begin{equation}\label{eq:porosity}
    \por(E)=\inf_{x\in E} \por(E,x).
\end{equation}
We say a set $E$ is \textit{porous} if $\por(E)>0$. Note that $0\leq \por(E) \leq \tfrac{1}{2}$, with $\por(E)$ close or equal to zero indicating there are few or no holes (e.g. $\por(\Q^{d})=0$) and $\por(E)$ close to $\tfrac{1}{2}$ indicating there are many large holes surrounding $E$ (e.g. for any hyperplane $\cH\subset \R^{d}$, $\por(\cH)=\tfrac{1}{2}$). As the examples illustrate a set can have relatively large packing dimension yet still be maximally porous, or have small packing dimension and be non-porous. Nevertheless, one can use porosity to obtain bounds on the packing dimension of sets, see \cite[Proposition 3.12]{Falconerbook}.

The main technical theorem, Theorem \ref{main}, of this article describes a general framework for studying the packing dimension of sets which avoid neighbourhoods of point collections of quantified size, by proving that sets of this type have positive {\it mean porosity}. Mean porosity generalises the above idea by allowing a set to be ``porous on average". Sets of this type arise naturally in Diophantine approximation in the study of badly approximable points, which we will now describe. 

The set of badly approximable points is a classical object of study in Diophantine approximation. Equipped with a norm $\|\cdot\|$ in $\R^{d}$, the set of \textit{$\varepsilon$-badly approximable points} is defined as
\begin{equation*}
    \Bad_{d}(\varepsilon,\|\cdot\|):=\left\{ \bx \in \R^{d} : \left\| \bx - \tfrac{\bp}{q} \right\| \geq \varepsilon q^{-1-\tfrac{1}{d}} \quad \text{ for all } (\bp,q)\in \Z^{d}\times \N  \right\},
\end{equation*}
and the set of \textit{badly approximable points} is defined as $$\Bad_{d}:=\bigcup_{\varepsilon>0} \Bad_{d}(\varepsilon,\|\cdot\|).$$ By the equivalence of norms in $\R^d$ note that $\Bad_{d}$ is independent of $\|\cdot\|$. It is well-known that $\Bad_{d}$, and by extension $\Bad_{d}(\varepsilon,\|\cdot\|)$, is a Lebesgue nullset, see for example \cite{BDGW}. To distinguish between the size of Lebesgue nullsets one can use various notions of dimension. In this article we are concerned with the packing dimension, denoted $\dimp$, and Hausdorff dimension, denoted $\dimh$. See \cite[\S 3]{Falconerbook} for definitions.

It has been shown that $\Bad_{d}$ is a set of full Hausdorff dimension. The case $d=1$ was proven by Jarnik \cite{Jarnik1928} using continued fraction theory, and for $d>1$ by Schmidt \cite{Schmidt1966,Schmidt1969}. Moreover, $\Bad_{d}$ is known to be \textit{hyperplane absolute winning}, see \cite{BFKRW12}. By properties of the packing dimension it immediately follows that $\dimp \Bad_{d} =d$.

The precise dimension of $\Bad_{d}(\varepsilon,\|\cdot\|)$ is far harder to calculate, so one usually is concerned with the asymptotic rate at which $\Bad_{d}(\varepsilon,\|\cdot\|)$ tends to full dimension as $\varepsilon\to 0$. Concerning the Hausdorff dimension, it was shown in \cite{Hensley1992,Simmons2018} that
\begin{equation} \label{eq: asymptotic dim of bad}
    \dimh \Bad_{d}(\varepsilon,\|\cdot\|) = d- V_{d,\|\cdot\|} \varepsilon^{d} + o(\varepsilon^{d}),
\end{equation}
for $V_{d,\|\cdot\|}$ an explicitly computable constant depending on the norm and dimension. For example, $V_{1}=\frac{6}{\pi^{2}}$. The case $d=1$ was first estimated by Kurzweil, before the sharp result \eqref{eq: asymptotic dim of bad} was given by Hensley \cite{Hensley1992} using continued fraction theory.
For $d>1$ estimates of $\dimh \Bad_{d}(\varepsilon,\|\cdot\|)$ appeared in \cite{SWeil2015,BroderickKleinbock2015}, before the sharp estimate \eqref{eq: asymptotic dim of bad} appeared in \cite{Simmons2018}.

In $d=1$ one can relate $\Bad_{1}(\varepsilon)$ to sets of real numbers with bounded partial quotients in their continued fraction expansion. It was shown in \cite{MaulUrbanPLMS, MaulUrbanTAMS} that the packing and Hausdorff dimensions of these sets are the same, and so the asymptotic rate of $\dimp \Bad_{1}(\varepsilon)$ is also known. In particular, for $\varepsilon>0$ sufficiently small
\begin{equation*}
    \codimp \Bad_{1}(\varepsilon)\asymp \codimh \Bad_{1}(\varepsilon).
\end{equation*}

 Under the general framework alluded to above we prove the following result.
\begin{theorem}\label{thm: packing dim bad}
    For $\varepsilon>0$ sufficiently small we have that 
    \begin{equation*}
        d-\dimp \Bad_{d}(\varepsilon, \|\cdot\|)\asymp \varepsilon^{d}.
    \end{equation*}
    In particular, for every $d\in \N$ and $\varepsilon>0$ sufficiently small
    \begin{equation*}
    \codimp \Bad_{d}(\varepsilon, \|\cdot\|) \asymp \codimh \Bad_{d}(\varepsilon, \|\cdot\|).
\end{equation*}
\end{theorem}

\begin{remark} 
The lower bound  on the dimension immediately follows from the \eqref{eq: asymptotic dim of bad}, the contribution of our general framework is the corresponding upper bound. It may be possible that upper bounds can be deduced from previous results. For example, while not stated, it can be seen from \cite{BroderickKleinbock2015} that
\begin{equation*}
   \dimp \Bad_{d}(\varepsilon,\|\cdot\|)\leq d - C \frac{\varepsilon^{d}}{|\log \varepsilon|},
\end{equation*}
for some constant $C>0$. It may be possible that the proof technique of \cite[Theorem 1.3]{Simmons2018} gives \eqref{eq: asymptotic dim of bad} for the packing dimension. The key contribution of this paper is the general Theorem~\ref{main} below.
\end{remark}

\begin{remark}
    The notion of porosity and badly approximable sets has appeared in the literature before, see 
    \cite{McMullen10, FRS16}. For instance, in \cite[Corollary 3.5]{FRS16}, it is shown that $\Bad_{d}(\varepsilon,\|\cdot\|)$ is a union of countably many $(\frac{\varepsilon}{3})^{d+1}$-porous sets. As can be seen in \cite[\S 3.8]{Falconerbook} porosity of this type gives the dimension bound
\begin{equation*}
    \dimp \Bad_{d}(\varepsilon, \|\cdot\|)\leq d-C_{d}\left(\frac{\varepsilon}{3}\right)^{d(d+1)}
\end{equation*}
for a constant $C_{d}>0$ depending on $d$ only.
\end{remark}

Using the general framework of Theorem \ref{main} alluded to above, in addition to Theorem~\ref{thm: packing dim bad}, we prove the following results in Diophantine approximation: 
\begin{enumerate}
    \item Upper bound to the asymptotic for the packing dimension of $\varepsilon$-Badly approximable points on manifolds (see Theorem \ref{thm:bad on manifold}).
    \item Upper bound to the asymptotic for the packing dimension of the twisted $\varepsilon$-badly approximable points (see Theorem \ref{bad irrational}). 
\end{enumerate}

The note is laid out as follows. Firstly we introduce our general framework on $\varepsilon$-badly approximable sets, and give our general result, see Theorem~\ref{main}. We then recall the notion of mean porosity before using this to prove Theorem~\ref{main}. In the final section we give some applications of Theorem~\ref{main}, including the proofs of Theorems~\ref{thm: packing dim bad}, \ref{thm:bad on manifold}, and \ref{bad irrational}. 

\subsection*{Acknowledgments} The authors would like to thank Sanju Velani for pointing out the work of \cite{FRS16} and Victor Beresnevich for helpful comments. The second name author is supported by a Leverhulme Early Career Research Fellowship (Grant no. ECF-2024-401). HK would like to acknowledge, that this material is based upon work supported by the Swedish Research Council under grant no. 2021-06594 while HK was in residence at Institut Mittag-Leffler in Djursholm, Sweden during the semester ``Interactions between fractal geometry, harmonic analysis, and dynamical systems" in autumn 2026. 

\subsection*{AI Declaration} AI has not been used in this work.





\section{Porosity and Badly approximable sets} \label{sec: main result} 

Consider the following setup: 
Let $\cR=(R_{\alpha})_{\alpha\in J}$ be a sequence of subsets $R_{\alpha}\subset \R^{d}$ indexed by an infinite countable set $J$. The set $R_{\alpha}$ are often called \textit{resonant sets}, see for example \cite{BDV06}. Let $h:J\to \R_{+}$ associate a height to each resonant set and assume that the height is unbounded on $J$ and that $\#\{\alpha \in J: h(\alpha)<N\}<\infty$ for every $N\in \N$. For $r>0$, denote the $r$-neighbourhood of a resonant set $R_{\alpha}$ as
\begin{equation*}
    \Delta(R_{\alpha},r):=\{\by \in \R^{d}: d(\by, R_{\alpha})<r\}, 
\end{equation*}
where for $\bx\in \R^d$ and $S\subset \R^d$
\begin{equation*}
    d(\bx, S):=\inf_{\by\in S} \|\bx-\by\|.
\end{equation*}
For an open set $\cU\subset \R^{d}$, and $\delta,\varepsilon>0$ define the set
\begin{equation*}
    \cB_{\cR}(\varepsilon)=\cB_{\cR}(\varepsilon,\delta,\cU):=\left\{ \bx \in \cU: \bx \not\in \Delta(R_{\alpha},\varepsilon h(\alpha)^{-\delta} ) \quad \text{ for all } \alpha \in J \right\}.
\end{equation*}
As an example of the framework, take $J=\N$ and $R_{q}=\left\{\frac{\bp}{q}: \bp \in \Z^{d}\right\}$, so $\cR=\Q^{d}$. Take $\delta=\frac{d+1}{d}$ and $h(q)=q$. Let $\cU=(0,1)^{d}$. Then we have that
\begin{equation*}
    \Bad_{d}(\varepsilon, \|\cdot\|)= \bigcup_{\bz \in \Z^{d}} \cB_{\Q^{d}}\left(\varepsilon, \tfrac{d+1}{d}, (0,1)^{d}+\bz\right).
\end{equation*}

The following is the abstract main result, from which all the applications follow. 

\begin{theorem} \label{main}
    Use the notation on resonance sets described above. Let $\cS \subseteq \N$ be an infinite set. Suppose that there exists constants $a\geq 0$, $b,c_{1},c_{2}>0$ such that for every $\bx \in \cU$ and every $\ell \in \cS$ there exists solution $\alpha \in J$ satisfying
    \begin{equation} \label{eq:dirichlet type cond}
            d(\bx, R_{\alpha})< \frac{c_{1}}{h(\alpha)^{a} 2^{\ell b}}, \qquad
            0\leq h(\alpha)< c_{2}2^{\ell}.
    \end{equation}
    Then, for $0<\varepsilon<c_{1}c_{2}^{b}$, we have that
    \begin{equation*}
     \dimp \cB_{\cR}(\varepsilon,a+b,\cU) \leq d- C \, \liminf_{M\to \infty}\frac{\#(\cS\cap [1,M])}{M}\left(\frac{\varepsilon}{c_{1}c_{2}^{b}}\right)^{d},
    \end{equation*}
    for constant $C>0$ dependent on $a,b,d$ only.
\end{theorem}

\begin{remark}
    By our assumption \eqref{eq:dirichlet type cond} it follows that for every $\bx \in \cU$ there exists $\alpha \in J$ such that
    \begin{equation*}
        d(\bx, R_{\alpha})< c_{1}c_{2}^{b} h(\alpha)^{-(a+b)},
    \end{equation*}
    and so $\cB_{\cR}(c_{1}c_{2}^{b},a+b,\cU)=\emptyset$. From the framework above there is no guarantee that there exists $\varepsilon>0$ such that $\cB_{\cR}(\varepsilon,a+b,\cU)\neq \emptyset$. In these instances the dimension bound would be trivial. In many of our applications there exists previously known results implying that $\cup_{\varepsilon>0}\cB_{\cR}(\varepsilon,a+b,\cU)$ has full Hausdorff dimension. In these instances the above dimension bound is non-trivial.
\end{remark}

    Other generalised notions of $\varepsilon$-badly approximable points, and how to calculate bounds on their Hausdorff dimension, exist. See work of \cite{Schmidt1966, KTV06, KleinbockWeiss2010, BFS19,SWeil2015,Tseng26} for various techniques to provide upper and lower bounds on the Hausdorff dimension of such sets. 
    
    

 The key step to proving Theorem~\ref{main} is to show that $\cB_{\cR}(\varepsilon,a+b,\cU)$ satisfy a version of \textit{porosity}. We give a brief introduction to packing dimension and mean porosity in \S~\ref{sec: packing and porosity} before proving Theorem~\ref{main}. We finish by demonstrating the breadth of our main theorem with a selection of applications, see \S~\ref{sec: applications}.

\subsection{Packing dimension and porosity} \label{sec: packing and porosity}

We begin by briefly recalling the definition of packing dimension. Introduced by Tricot \cite{Tricot82}, for $E\subseteq \R^{d}$ the packing dimension of $E$ is defined to be
\begin{equation*}
    \dimp E := \inf \left\{ \sup_{i} \left(\limsup_{r\to 0} \frac{\log_{2} N_{r}(E_{i})}{-\log_{2} r} \right) : E=\bigcup_{i} E_{i} \right\},
\end{equation*}
where $N_{r}(E_{i})$ is the smallest number of cubes of sidelength $r$ required to cover $E_{i}$. Note one can equivalently define $\dimp E$ in terms of packing measures (defined similarly to Hausdorff measures, but one takes packings of a set rather than covers), see \cite[\S 3.5]{Falconerbook} for details. Packing dimension enjoys many of the properties that Hausdorff dimension possesses. In particular, packing dimension is;
\begin{itemize}
    \item  \textit{monotonic:} If $E\subseteq F$ then $\dimp E \leq \dimp F$.
    \item \textit{Countably stable:} For all countable collections $\{F_{i}\}_{i}$ of subsets of $\R^{d}$ we have $\dimp \bigcup_{i} F_{i} = \sup_{i} \dimp F_{i}$.
    \item \textit{Bi-Lipschitz stable:} For subset $F\subset \R^{d}$, if $f:F \to \R^{m}$ is a bi-Lipschitz transformation, that is
    \begin{equation*}
        c_{1}|\bx-\by| \leq |f(\bx)-f(\by)|\leq c_{2}|\bx-\by| \quad \text{ for all }\, \bx,\by \in F,
    \end{equation*}
    then $\dimp f(F)=\dimp F$.
\end{itemize}
In comparison to Hausdorff dimension we have the useful inequality that for any set $F\subset \R^{d}$
\begin{equation*}
    \dimh F \leq \dimp F.
\end{equation*}
 See for example \cite[\S 2.4]{Fraserbook} for a list of various properties between dimensions. There are many examples of sets where we have equality, and many where the inequality is strict. 

Recall the definition of porosity from \eqref{eq:porosity}. We now formally define the notion of mean porosity. Let $\alpha \in [0,\tfrac{1}{2}]$ and $\rho\in [0,1]$. We say a set $E\subset \R^{d}$ is $(\alpha,\rho)$-mean porous if, for all $x\in E$,
\begin{equation*}
    \liminf_{N\to\infty} \frac{1}{N}\# \left\{ 1\leq i \leq N: \por(E,x,2^{-i})\geq \alpha \right\} \geq \rho.
\end{equation*}
We say a set $E$ is \textit{mean porous} if there exists a pair $(\alpha,\rho)\in (0,\tfrac{1}{2}]\times (0,1]$ such that $E$ is $(\alpha,\rho)$-mean porous. Mean porous sets have the following packing dimension bound.

\begin{theorem}\label{mean porous dim}
    Suppose $E\subset \R^{d}$ is $(\alpha,\rho)$-mean porous. Then
    \begin{align*}
        \dimp E &\leq d-C_{d}\rho  \alpha^{d}, \\
        \dimp E &\leq d-\rho + \frac{C_{d}}{|\log(1-2\alpha)|},
    \end{align*}
    for constant $C_{d}>0$ dependent only on $d$.
\end{theorem}

The version cited above appears in \cite[Theorem 8]{ShmerkinSurvey}. The first bound, due to Koskela and Rohde \cite{KosRoh97}, is generally good for porosity close to zero and the second bound, see for example \cite{BJJKRSS09}, is generally good for porosity close to $\tfrac{1}{2}$. Note that the second bound gives a packing dimension upper bound of $d-1$ on maximally porous sets (seen to be $(\tfrac{1}{2},1)$-mean porous). 

We remark that there are many other notions of porosity, such as porosity of measures and porosity for general metric spaces \cite{JJKRRS10, ShmerkinSurvey}, although those won't be discussed here. 

\section{Proof of Theorem~\ref{main}}
Following from the previous section, in order to prove Theorem~\ref{main} we show the following result.

\begin{theorem}\label{thm: porous statement}
    Let $\cS \subseteq \N$ be an infinite set. Suppose that there exists constants $a\geq 0$, $b,c_{1},c_{2}>0$ such that for every $\bx \in \cU$ and every $\ell \in \cS$ there exists $\alpha \in J$ such that
    \begin{equation} \label{eq:dirichlet type cond thm}
            d(\bx, R_{\alpha})< \frac{c_{1}}{h(\alpha)^{a} 2^{\ell b}}, \qquad
            0\leq h(\alpha)< c_{2}2^{\ell}.
    \end{equation}
    Then, for every $0<\varepsilon<c_{1}c_{2}^{b}$, the set $\cB_{\cR}(\varepsilon,a+b,\cU)$ is 
    \begin{equation*}
        \left( \frac{\varepsilon}{8c_{1}c_{2}^{b}}, \, \frac{1}{8\lceil b^{-1} \rceil (a+b)}\liminf_{M\to \infty}\frac{\#(\cS\cap [1,M])}{M} \right)\text{-mean porous.}
    \end{equation*}
\end{theorem}

By combining Theorem~\ref{thm: porous statement} with Theorem~\ref{mean porous dim} we immediately obtain that
\begin{equation*}
        \dimp \cB_{\cR}(\varepsilon,a+b) \leq d- \frac{C_{d}}{8^{d+1}\lceil b^{-1}\rceil (a+b)}\liminf_{M\to \infty}\frac{\#(\cS\cap [1,M])}{M}\left(\frac{\varepsilon}{c_{1}c_{2}^{b}}\right)^{d},
\end{equation*}
thus completing the proof of Theorem~\ref{main} as a consequence of Theorem~\ref{thm: porous statement}.

\begin{proof}[Proof of Theorem~\ref{thm: porous statement}]
For ease of notation set
\begin{equation*}
    \cB_{\cR}(\varepsilon):=\left\{ \bx \in \cU : \bx \not \in \Delta\left(R_{\alpha}, \varepsilon h(\alpha)^{-a-b}\right) \quad \text{ for all} \, \, \alpha \in J \right\}.
\end{equation*}
Set
\begin{equation*} 
    D=2\left\lceil \frac{1}{b} \right\rceil.
\end{equation*}
    Let $\widehat{\cS}$ be a maximal $D$-separated subset of $\cS$ in the sense that for every $l,l'\in\widehat{\cS}$ we have 
    \begin{equation*}
    |l-l'|\geq D,    
    \end{equation*}
    and for every $l\in \cS$ there exists $l'\in \widehat{\cS}$ such that
    \begin{equation*}
        |l-l'|< D.
    \end{equation*}    
    If $\cB_{\cR}(\varepsilon)=\emptyset$ the dimension bound is trivial, so assume $\cB_{\cR}(\varepsilon)$ is non-empty. Take any $\bx \in \cB_{\cR}(\varepsilon)$ and any $r \in \widehat{\cS}$. Applying statement \eqref{eq:dirichlet type cond thm} we see that there exists solution $\alpha\in J$ and resonant set $R_{\alpha} \in \cR$ solving
    \begin{equation} \label{eq: x-r ball}
        d(\bx, R_{\alpha})<c_{1} h(\alpha)^{-a}2^{-rb} , \qquad 0 \leq h(\alpha)<c_{2} 2^{r}.
    \end{equation}
    Suppose $R_{\alpha(r)}$ is the resonant set with smallest height that satisfies \eqref{eq: x-r ball}. That is,
    \begin{equation*}
        h(\alpha(r)):= \min \left\{ h(\alpha) : \alpha \in J \, \text{ satisfying \eqref{eq: x-r ball} } \right\}.
    \end{equation*}
    The minimum is well defined by the fact that $1\leq \#\{\alpha \in J: h(\alpha)<c_{2}2^{r}\}<\infty$. 
    Let $k_{r}\in \Z$ be the unique integer such that
    \begin{equation} \label{eq:k_r bound}
        2^{-k_{r}-1}\leq c_{1} h(\alpha(r))^{-a} < 2^{-k_{r}}.
    \end{equation}
    Note by the minimality of the height of $\alpha(r)$ satisfying \eqref{eq: x-r ball} that $k_{r}$ is non-decreasing as $r$ increases. To see why suppose that $a>0$ and $r'>r$. Then comparing $\alpha(r)$ and $\alpha(r')$ we see that if $h(\alpha(r'))<h(\alpha(r))$, then
    \begin{equation*}
        d(\bx, R_{\alpha(r')})<c_{1}h(\alpha(r'))^{-a}2^{-r'b}<c_{1}h(\alpha(r'))^{-a}2^{-rb}, \qquad h(\alpha(r'))<h(\alpha(r))<c_{2}2^{r}
    \end{equation*}
    meaning $\alpha(r')$ satisfies \eqref{eq: x-r ball} thus contradicting the minimality of $h(\alpha(r))$. If $a=0$, then $k_{r}\in \Z$ is a fixed integer independent of $r$. Thus, for all $r'>r$ we have $h(\alpha(r'))\geq h(\alpha(r))$. The non-decreasing property of $k_{r}$ immediately follows.
    
    Using \eqref{eq: x-r ball} and \eqref{eq:k_r bound} we can see that
    \begin{equation} \label{eq: height of k_r bound}
        k_{r}<a\log_{2}h(\alpha(r))-\log_2 c_1 < ra +a\log_{2}c_{2}-\log_{2}c_{1}.
    \end{equation}
    It follows from \eqref{eq: x-r ball} and \eqref{eq:k_r bound} that there exists some point $\by(r) \in R_{\alpha(r)}$ such that $\by(r) \in B(\bx, 2^{-k_{r}-rb})$. By definition of $\cB_{\cR}(\varepsilon)$, we have that 
    \begin{equation*}
        \cB_{\cR}(\varepsilon)\cap \Delta\left(R_{\alpha(r)}, \varepsilon h(\alpha(r))^{-a-b}\right)=\emptyset,
    \end{equation*}
    and since $\by(r)\in R_{\alpha(r)}$ it follows that
\begin{equation*}
    B(\by(r),\varepsilon h(\alpha(r))^{-a-b}) \subset \Delta(R_{\alpha(r)}, \varepsilon h(\alpha(r))^{-a-b})
\end{equation*}
    and so
    \begin{equation*}
    B(\by(r), \varepsilon h(\alpha(r))^{-a-b})\cap \cB_{\cR}(\varepsilon)=\emptyset.
    \end{equation*}
    Therefore, there exists some $\bz\in B(\by(r), \varepsilon h(\alpha(r))^{-a-b}))$ such that
    \begin{equation*}
        B\left(\bz, \tfrac{1}{2}\varepsilon h(\alpha(r))^{-a-b}\right) \subset B(\bx, 2^{-k_{r}-rb})\setminus \cB_{\cR}(\varepsilon).
    \end{equation*}
    Taking the slightly larger ball $B(\bx, 2^{-\lfloor k_{r}+rb\rfloor})$, we have that
    \begin{equation} \label{eq: por scales}
        \por\left(\cB_{\cR}(\varepsilon),\bx, 2^{-\left\lfloor k_{r}+rb\right\rfloor }\right) \geq \frac{\frac{1}{2}\varepsilon h(\alpha(r))^{-a-b}}{2^{-k_{r}-\lfloor rb \rfloor }} \overset{\text{\eqref{eq:k_r bound}}}{\geq} \frac{1}{4c_{1}}\varepsilon\left(\frac{2^{\frac{1}{b}\lfloor rb \rfloor}}{h(\alpha(r))}\right)^{b}    \overset{\text{\eqref{eq: x-r ball}}}{>} \frac{\varepsilon }{8c_{1}c_{2}^{b}}.
    \end{equation}
    This holds for every $r\in \widehat{\cS}$. For ease of notation let 
    \begin{equation*}
        t(r)=\lfloor k_{r}+rb\rfloor, \qquad r \in \widehat{\cS}.
    \end{equation*}
    These are the dyadic scales at which we find suitably large holes. Note the following two properties of the sequence $(t(r))_{r\in \widehat{\cS}}$:
    \begin{itemize}
        \item[(P1)] \textit{ For any $r,r'\in \widehat{\cS}$ with $r\neq r'$ we have that $t(r)\neq t(r')$}: for any two $r,r' \in \widehat{\cS}$ with $r> r'$ we have that $t(r)> t(r')$. To see why observe that if $a=0$ then $k_{r}=k_{r'}$, otherwise use that $k_{r}$ is non-decreasing, so as $r>r'$ it follows that $k_{r}\geq k_{r'}$. By the $D$-separated property of $\widehat{\cS}$ we have $\lfloor rb \rfloor -\lfloor r'b \rfloor >0$. Combining we obtain that $\lfloor k_{r}+rb \rfloor> \lfloor k_{r'}+r'b \rfloor$ as required.\\

        \item[(P2)] \textit{$t(r)$ is bounded linearly from above by $r$}: Using \eqref{eq: height of k_r bound} we have that
        \begin{equation*}
            t(r)\leq k_{r}+rb <r(a+b)+a\log_{2}c_{2}-\log_{2}c_{1}.
        \end{equation*}
    \end{itemize}
    
    We now calculate the frequency of the dyadic scales at which we can find these holes. It follows from the maximal property of $\widehat{\cS}$ in $\cS$ that for any sufficiently large $R\in \R_{>0}$ we have
\begin{align} \label{eq: bound on S hat}
        \#\left(\cS\cap \left[1,\tfrac{1}{2}R\right] \right) &\leq \#\left(\cS\cap \left[1,R-D\right] \right) \nonumber \\
        &\leq 2D\#\left(\widehat{\cS}\cap \left[1,R\right]\right).
\end{align}

    For any large $M\in \N$ we have that
    \begin{align*}
        \frac{1}{M}\#\left\{ 1\leq n\leq M: \por(\cB(\varepsilon),\bx,2^{-n})\geq \tfrac{1}{8c_{1}c_{2}^{b}}\varepsilon \right\} &\overset{\text{(P1)} + \eqref{eq: por scales}}{\geq} \frac{1}{M} \#\{r \in \widehat{\cS} : 1\leq t(r)\leq M\}\\
        & \hspace{-5cm}\overset{\text{(P2)}}{\geq} \frac{1}{M}\left(\#\left(\widehat{\cS}\cap \left[ 1, \frac{M-a\log_{2}c_{2}+\log_{2}c_{1}}{a+b}\right] \right)-\#\left\{r \in \widehat{\cS}: t(r)<1\right\}\right)\\
        &\hspace{-5cm}\overset{\text{\eqref{eq: bound on S hat}}}{\geq} \frac{1}{M}\left(\frac{1}{2D}
        \#\left(\cS\cap \left[1, \frac{M-a\log_{2}c_{2}+\log_{2}c_{1}}{2(a+b)} \right]\right)-\#\left\{r \in \widehat{\cS}: t(r)<1\right\}\right)
    \end{align*}
    Since $c_{2}, c_{1}, a$ and $\#\{r \in \widehat{\cS}:t(r)<1\}<\infty$ (since $t(r)\to\infty$ and $r\to \infty$) is fixed independent of $M$, we have that
    \begin{align*}
        \liminf_{M\to\infty}& \frac{1}{M}\#\left\{1\leq n\leq M: \por(\cB(\varepsilon),\bx,2^{-n})\geq \tfrac{1}{8c_{1}c_{2}^{b}}\varepsilon\right\} 
        \geq \frac{1}{4D(a+b)}\liminf_{M\to\infty}\frac{\#(\cS\cap [1,M])}{M}.
    \end{align*}
    This is true for any $\bx \in \cB(\varepsilon)$ and so $\cB(\varepsilon)$ is $\left( C_{1}\varepsilon, C_{2}\liminf_{M\to\infty}\frac{\#(\cS\cap [1,M])}{M}\right)$-mean porous, for 
    \begin{align*}
        C_{1}&=\frac{1}{8c_{1}c_{2}^{b}},\\
        C_{2}&=\frac{1}{8\lceil b^{-1}\rceil (a+b)},
    \end{align*}
    as required.
\end{proof}

\section{Applications} \label{sec: applications}
We now apply our main results to several different problems. For completeness we begin with the easy case of proving Theorem~\ref{thm: packing dim bad}

\subsection{$\varepsilon$-Badly approximable points}
Note that the lower bound on the packing dimension can be deduced from the result of Simmons \cite[Theorem 1.3]{Simmons2018}. In particular, for $\varepsilon>0$ suffciently small one has that
\begin{equation*}
    \dimh \Bad_{d}(\varepsilon, \|\cdot\|)\geq d- 2 V_{d,\|\cdot\|}\varepsilon^{d}.
\end{equation*}
Since $\dimh \Bad_{d}(\varepsilon, \|\cdot\|)\leq \dimp \Bad_{d}(\varepsilon, \|\cdot\|)$ we have that
\begin{equation*}
    d-\dimp \Bad_{d}(\varepsilon, \|\cdot\|)\leq 2V_{d,\|\cdot\|} \varepsilon^{d}.
\end{equation*}
Hence Theorem~\ref{thm: packing dim bad} follows from
\begin{theorem}
    For any $\varepsilon>0$ sufficiently small
    \begin{equation*}
        \dimp \Bad_{d}(\varepsilon, \|\cdot\|)\leq d-\frac{C_{d}}{(d+1)8^{d+1}c_{\|\cdot\|}^{d}} \varepsilon^{d}.
    \end{equation*}
\end{theorem}
\begin{proof}
Take
\begin{equation*}
     J=\N, \quad R_{q}=\left\{ \frac{\bp}{q}: \bp \in \Z^{d}\right\}, \quad h(R_{q})=q, \\
\end{equation*}
so $\cR=\Q^{d}$. By Minkowski convex body Theorem for any norm $\|\cdot\|$ there exists a constant $c_{\|\cdot\|}>0$ such that for every $\bx\in \R^{d}$ and every $\ell \in \N$ there exists solution $\frac{\bp}{q}\in \Q^{d}$ satisfying
\begin{equation*}
    \left\|\bx - \frac{\bp}{q}\right\|< \frac{c_{\|\cdot\|}}{q 2^{\ell \tfrac{1}{d}}}, \quad 1 \leq q < 2^{\ell}.
\end{equation*}
Take $\cU=(0,1)^{d}+\bz$ for any $\bz \in \Z^{d}$. Then \eqref{eq:dirichlet type cond thm} is satisfied with 
\begin{equation*}
    a=1, \quad b=\frac{1}{d}, \quad c_{1}=c_{\|\cdot\|}, \quad c_{2}=1,
\end{equation*}
and $\cS=\N$. Applying Theorem~\ref{main} immediately implies that for any $0<\varepsilon<c_{\|\cdot\|}$
\begin{equation*}
    \dimp \cB_{\Q^{d}}(\varepsilon, \tfrac{d+1}{d}, (0,1)^{d}+\bz)\leq d-\frac{C_{d}}{(d+1)8^{d+1}c_{\|\cdot\|}^{d}} \varepsilon^{d}.
\end{equation*}
Note that here we have kept explicit the dependence on $a,b$ so that $C_d$ is the constant from Theorem \ref{mean porous dim}. Taking the union over all $\bz \in \Z^{d}$ and using the countable stability of the packing dimension we obtain the desired result.
\end{proof}

\subsection{Badly approximable points on manifolds}

Motivated by conjecture's of Mahler (1932) and Sprindzhuk (1980), in the past thirty years there has been a significant amount of research in Diophantine approximation on manifolds, see \cite{kleinbockMargulis98, Beresnevich15} for key breakthrough results. Within the remit of this paper we are concerned with set of $\varepsilon$-badly approximable points intersected with a submanifold $\cM$ of $\R^{d}$.

It was shown by Baker in the 1970s that for any $C^{1}$ submanifold $\cM$ of $\R^{d}$ the set $\Bad_{d}\cap \cM$ is a nullset, with respect to the induced $\dimh \cM$-dimensional Lebesgue measure on $\cM$ \cite{Baker76}. Therefore, to give any meaningful notion of size the dimension of the set should be studied. For $d=2$ it was proven in \cite{BadVel14} that for any $C^{2}$ planar curve $\cC_{f}:=\{(x,f(x)): x\in I\}$, with $f:I\to\R$ a twice continuously differentiable function defined on an interval $I$, that is not a line segment
\begin{equation*}
    \dimh \Bad_{2}\cap \cC_{f} =1.
\end{equation*}
Note that the removal of the case $\cC_{f}$ a line segment is necessary. For example the statement is false for the straight line $\cL=\{(x,0): x\in I\}$, since $\Bad_{2}\cap \cL=\emptyset$. In higher dimensions this curvature condition is replaced with \textit{non-degeneracy}, see \cite{kleinbockMargulis98}. Imposing the additional condition of analyticity Beresnevich was able to prove the higher dimensional analogue \cite{Beresnevich15}. That is, for $\cM$ a non-degenerate analytic manifold
\begin{equation*}
    \cM:=\left\{ (\bx,f(\bx)): \bx \in \cU \subset \R^{m}\right\}\subset \R^{d}
\end{equation*}
with $\dimh \cM=m$ and codimension $n=d-m$, described by an analytic map $f=(f_{1},...,f_{n}):\cU\to \R^{n}$ Beresnevich proved that
\begin{equation*}
    \dimh \Bad_{d}\cap \cM = \dimh \cM=m.
\end{equation*}
Note this immediately implies $\dimp \Bad_{d}\cap \cM =m$. See \cite{Yang19,DattaShao26,Zhao26} for more recent results on $\Bad_{d}\cap \cM$, including winning properties, weighted/inhomogeneous analogues, and weakening on various conditions of $\cM$. 

A natural follow up question to the results of \cite{Beresnevich15,BadVel14} is to ask whether we can say anything about the smaller set $\Bad_{d}(\varepsilon, \|\cdot\|)\cap \cM$. We prove the following result.
\begin{theorem} \label{thm:bad on manifold}
Let $\cM=\{(\bx,f(\bx)): \bx \in \cU\subset \R^{m}\}\subset \R^{d}$ for $\cU$ an open ball with $f=(f_{1},\ldots, f_{n})$ twice continuously differentiable. Suppose that the first and second order partial derivatives of $f$ are bounded on $\cU$. That is,  
\begin{align*}
    F:=\underset{1\leq i \leq m}{\max_{1\leq j \leq n}} \sup_{\bx \in \cU} \left| \frac{\partial f_{j}}{\partial x_{i}}(\bx)\right|<\infty,\qquad \text{ and } \qquad
    S:= \underset{1\leq i,k \leq m}{\max_{1\leq j \leq n}} \sup_{\bx \in \cU} \left| \frac{\partial^{2} f_{j}}{\partial x_{i}\partial x_{k}}(\bx)\right|<\infty.
\end{align*}
Then, for $\varepsilon>0$ sufficiently small we have
\begin{equation*}
    \dimp \Bad_{d}(\varepsilon, \|\cdot\|_{\infty})\cap \cM \leq \dimh \cM -C_{\cM} \varepsilon^{d},
\end{equation*}
with constant $C_{\cM}>0$ depending on $d$ and $\cM$ only.
\end{theorem}

\begin{remark}\rm
 As far as the authors are aware prior to our result it was unknown whether it was possible to construct a submanifold $\cM$ such that $\dimh \Bad_{d}(\varepsilon, \|\cdot\|)\cap \cM=\dimh \cM$. Our result shows this is not possible in the case that $\cM$ is twice continuous differentiable with bounded first and second partial derivatives. Note that the non-degeneracy condition is not present. This is because we are only concerned with the dimension upper bound. Clearly, as the example $\cL$ above illustrates, there are instances of degenerate manifolds where our bound would be trivial. 
\end{remark}

Observe that the exponent on $\varepsilon>0$ is $d$, independent of $\dimh \cM$. One may ask whether this is the optimal bound. By combining the result of Simmons on $\dimh \Bad_{d}(\varepsilon,\|\cdot\|)$, \cite[Theorem 1.3]{Simmons2018} with \cite[Theorem 8.2]{Falconerbook} one can deduced the following easy result.

\begin{proposition} \label{prop: bad on typical manifold}
    Fix an $\varepsilon>0$ small and let $\cM$ be any manifold. Let $G$ denote the group of rigid motions on $\R^{d}$. Then for a positive measure set of $f\in G$ we have
    \begin{equation*}
        \dimh \Bad_{d}(\varepsilon, \|\cdot\|_{\infty})\cap f(\cM) \geq \dimh \cM-V_{d} \varepsilon^{d} +o(\varepsilon^{d}).
    \end{equation*}
\end{proposition}
That is, for every $\varepsilon>0$ sufficiently small there exists manifolds $\cM$ for which the exponent of $d$ on the $\varepsilon$ term in Theorem \ref{thm:bad on manifold} is optimal. We make no claim on the optimality of the multiplicative constant on $\varepsilon^{d}$ appearing in Theorem~\ref{thm:bad on manifold}.
\subsubsection{Proof of Theorem~\ref{thm:bad on manifold}}
We now turn to the proof of Theorem~\ref{thm:bad on manifold}. We need the following preliminary result that can be found in \cite[Theorem 4]{BLVV17}. It essentially gives us a Dirichlet-type theorem on manifolds.
\begin{theorem}[{\cite[Theorem 4]{BLVV17}}] \label{thm:BLVV}
    Let $\cM=\{(\bx,f(\bx)): \bx \in \cU\subset \R^{m}\}\subset \R^{d}$ for $\cU$ an open ball and $f=(f_{1},\ldots, f_{n}):\cU \to \R^{n}$ twice continuously differentiable with uniformly bounded first and second order partial derivatives. Fix $\eta>0$ and some ball $B\subset \cU$ with $\dist(B,\cU^{c})>0$. Then, for any $N>N_{0}=N_{0}(\eta,\cM,B)$ and any $\bx \in B$ there exists $(\bp_{1},\bp_{2},q)\in \Z^{m}\times \Z^{n}\times \N$ with $\tfrac{\bp_{1}}{q}\in \cU$ and $1 \leq q \leq N$ satisfying
    \begin{align*}
    \begin{cases}
        &\left\| \bx - \frac{\bp_{1}}{q}\right\|_{\infty} < \left( \frac{2}{\eta}\right)^{\tfrac{n}{m}} q^{-1} N^{-\tfrac{1}{d}}, \\
        &\left\| f\left(\frac{\bp_{1}}{q}\right)-\frac{\bp_{2}}{q} \right\|_{\infty}< \eta q^{-1-\tfrac{1}{d}}.
        \end{cases}
    \end{align*}
\end{theorem}
\begin{remark}
    In \cite{BLVV17} the theorem is applicable to general functions $\psi:\N \to \R_{+}$ satisfying
    \begin{equation*}
        \inf_{Q\in \N} Q\psi(Q)>0 \qquad \limsup_{Q\to \infty} Q^{1/n}\psi(q)=\infty.
    \end{equation*}
    The bound on $N$ appearing in \cite[Theorem 4]{BLVV17} is dependent on the choice of $\bx$. However, carefully going through the proof one sees that $N$ depends on $\dist(\bx,\cU^{c})$. Precisely, for a fixed $\bx \in \cU$ take any ball $B(\bx,r)\subset \cU$ with $r>0$. Then the statement holds for all
    \begin{equation*}
        N\in \left\{ N\in \N: (2^{-n}N\psi(N)^{n})^{-\frac{1}{d}}<\min \left\{ 1,r, \left(\frac{\inf_{Q\in \N} Q\psi(Q)}{Sm^2}\right)^{1/2} \right\} \right\},
    \end{equation*}
    where $S$ is as in Theorem~\ref{thm:bad on manifold}. Since we are applying the result to functions of the form $\psi(Q)=\eta Q^{-\frac{1}{d}}$, and $n<d$, these conditions on the approximation function are satisfied. Moreover, the above set contains all integers $N \in \N$ such that
    \begin{equation*}
        N> N_{0}(\eta, \cM, \bx):=\left(\left(\frac{\eta}{2} \right)^{\frac{m}{n}}\min \left\{ 1,r, \left(\frac{\eta}{Sm^2}\right)^{1/2} \right\} \right)^{-\frac{md}{d-n}}.
    \end{equation*}
    To remove the dependency on $\bx$ note that for any $\bx \in B\subset \cU$ with $\dist(B, \cU^{c})=r_{B}>0$ it follows that $B(\bx,r_{B})\subset \cU$, so the choice of $r$ can be made dependent on $B$. Thus the Theorem holds for all $N\geq N_{0}$ with
    \begin{equation*}
        N_{0}=N_{0}(\eta, \cM, B):=\left(\left(\frac{\eta}{2} \right)^{\frac{m}{n}}\min \left\{ 1,r_{B}, \left(\frac{\eta}{Sm^2}\right)^{1/2} \right\} \right)^{-\frac{md}{d-n}}.
    \end{equation*}
\end{remark}

\begin{proof}[Proof of Theorem~\ref{thm:bad on manifold}]
Fix $\varepsilon>0$. Choose parameters $1>\rho,\delta>0$ such that
\begin{equation} \label{eq: rho delta choice}
    Fm\rho+\delta<1.
\end{equation}
Note these are independent of $\varepsilon$. Let $(B_{i})_{i\in \N}$ be a countable collection of open balls such that $\overline{B_{i}} \subset \cU$ and 
\begin{equation*}
    \bigcup_{i}B_{i}=\cU.
\end{equation*}
Then $\dist(B_{i},\cU^{c})=r_{i}>0$ for all $i\in \N$. Fix some $B_{i}$ and define the set 
\begin{equation*}
    \cR(\delta \varepsilon):=\left\{ \frac{\bp_{1}}{q}\in \Q^{m}: \frac{\bp_{1}}{q}\in \cU \quad \text{ and } \quad \left\| f\left(\frac{\bp_{1}}{q}\right)-\frac{\bp_{2}}{q} \right\|< \delta \varepsilon q^{-1-\tfrac{1}{d}} \quad \text{ for some } \bp_{2}\in \Z^{n} \right\}\subset \R^{m}\cap \cU
\end{equation*}
As in the notation of \S~\ref{sec: main result} set 
\begin{align*}
    J= \N,\qquad R_{q}=\left\{\frac{\bp}{q}\in \cR(\delta\varepsilon)\right\}, \qquad h(R_{q})=q.
\end{align*}
In this application we may take $\cS=\N_{\geq \log_{2}N_{0}}$ with $N_{0}$ depending on $\delta\varepsilon$, $\cM$ and $B_{i}$.
Theorem~\ref{thm:BLVV} with $\eta=\delta \varepsilon$ implies condition \eqref{eq:dirichlet type cond} is satisfied with parameters
\begin{align*}
    a=1, \qquad b=\frac{1}{d}, \qquad c_{1}=\left(\frac{2}{\delta \varepsilon}\right)^{\frac{n}{m}}, \qquad c_{2}=1.
\end{align*}
Note that $c_{1}$ is dependent on $\varepsilon>0$. Applying Theorem~\ref{main} to the set $\cB_{\cR(\delta\varepsilon)}\left(\rho\varepsilon,\tfrac{d+1}{d},B_{i}\right)$ we have, for any $0<\varepsilon<\left(\frac{2}{\delta \varepsilon}\right)^{\frac{n}{m}}$ or equivalently
\begin{equation*}
    0<\varepsilon<\left(\frac{2}{\delta}\right)^{\frac{n}{d}},
\end{equation*}
the dimension bound
\begin{align} \label{eq: dim result step 1}
    \dimp \cB_{\cR}\left(\rho \varepsilon, \tfrac{d+1}{d},B_{i}\right) &\leq m- C \left(\frac{\rho\varepsilon}{c_{1}}\right)^{m}\nonumber\\
    &\leq m-C\left(\frac{\delta}{2}\right)^{n}\rho^{m}\varepsilon^{m+n}\nonumber \\
    &\leq m-C_{\cM} \varepsilon^{d}, \qquad 
\end{align}
with $C_{\cM}=C \rho^{m}\left(\frac{\delta}{2}\right)^{n}$. Note that $\rho$ and $\delta$ depend on $F$ and thus $\cM$ and importantly are independent of $B_{i}$. Moreover,
\begin{equation*}
    \liminf_{M\to\infty}\frac{M-\log_{2}N_{0}}{M}=1,
\end{equation*}
independent of $\delta\varepsilon$, $\cM$ and $B_{i}$.

Let $\pi_{m}:\R^{d}\to \R^{m}$ denote the orthogonal projection of $\R^{d}$ onto the first $m$ coordinates. Let
\begin{equation*}
    \cM_{i}:=\{(\bx,f(\bx)): \bx \in B_{i}\}
\end{equation*}
Assume that
\begin{equation} \label{eq: bi-lip property}
    \pi_{m}\left(\Bad_{d}(\varepsilon,\|\cdot\|_{\infty})\cap \cM_{i}\right) \subseteq \cB_{\cR(\delta\varepsilon)}\left(\rho \varepsilon, \tfrac{d+1}{d},B_{i}\right).
\end{equation}
The map $\pi_{m}:\cM_{i} \to B_{i}$ is a bi-Lipschitz transformation since the partial derivatives are bounded on $\cU\supset B_{i}$. Since $\dimp$ is preserved under bi-Lipschitz maps and is countably stable we have that
\begin{align*}
    \dimp \Bad_{d}(\varepsilon, \|\cdot\|_{\infty})\cap \cM &=\sup_{i \in \N}\dimp \pi_{m}\left(\Bad_{d}(\varepsilon, \|\cdot\|_{\infty})\cap \cM_{i}\right) \\
    &\overset{\text{\eqref{eq: bi-lip property}}}{\leq} \sup_{i\in \N}\dimp \cB_{\cR(\delta\varepsilon)}\left(\rho\varepsilon, \tfrac{d+1}{d},B_{i}\right)\\
    &\overset{\text{\eqref{eq: dim result step 1}}}{\leq} \dimh \cM - C_{\cM} \varepsilon^{d}.
\end{align*}
Thus it remains to show \eqref{eq: bi-lip property}. Suppose that $\bx \in B_{i}\setminus \cB_{\cR(\delta\varepsilon)}(\rho \varepsilon,\tfrac{d+1}{d},B_{i})$. Then, there exists $\frac{\bp_{1}}{q}\in \cR(\delta\varepsilon)$ such that
\begin{equation*}
    \left\|\bx-\frac{\bp_{1}}{q} \right\|< \rho \varepsilon q^{-1-\tfrac{1}{d}}.
\end{equation*}
By definition of $\cR(\delta\varepsilon)$ there exists $\bp_{2}\in \Z^{n}$ such that
\begin{equation*}
    \left\| f\left( \frac{\bp_{1}}{q}\right) - \frac{\bp_{2}}{q} \right\|_{\infty} < \delta \varepsilon q^{-1-\tfrac{1}{d}}.
\end{equation*}
Now, by the triangle inequality, the Mean Value Theorem, and the fact that the first partial derivatives of $f$ are bounded from above by $F$ in $\cU$, we have that
\begin{align*}
    \left\| f(\bx)-\frac{\bp_{2}}{q} \right\|_{\infty} &\leq \left\| f(\bx)-f\left(\frac{\bp_{1}}{q}\right)\right\|_{\infty}+\left\| f\left(\frac{\bp_{1}}{q}\right)-\frac{\bp_{2}}{q}\right\|_{\infty}\\
    & < Fm\left\| \bx - \frac{\bp_{1}}{q} \right\|_{\infty} + \delta \varepsilon q^{-1-\tfrac{1}{d}}\\
    &\leq Fm \rho \varepsilon q^{-1-\tfrac{1}{d}} + \delta \varepsilon q^{-1-\tfrac{1}{d}} \\
    &< \varepsilon q^{-1-\tfrac{1}{d}}
\end{align*}
where the last line follows from \eqref{eq: rho delta choice}.
This implies that $\bx \not\in \pi_{m}\left(\Bad_{d}(\varepsilon, \|\cdot\|_{\infty})\cap \cM_{i}\right)$, and so \eqref{eq: bi-lip property}, and hence the theorem, is proven.
\end{proof}

\subsection{Twisted $\varepsilon$-Badly approximable points}
For a dynamical system $(T,X)$ the study of the \textit{shrinking target problem} looks at the set of points that are visited by shrinking neighbourhoods of the orbit $(T^{n}(\bx))_{n \in \N}$ of some fixed point $\bx \in X$. A classical example of such system is the orbit of irrational rotations on the unit circle as studied by Kurzweil \cite{Kurzweil1951} and since developed in numerous directions \cite{Kur55, Bug03,TKim24,HussainWard24}. In this application we study a natural higher dimensional analogue of irrational rotations.

Fix some $A\in \R^{m\times n}$ 
and $\varepsilon>0$. Define the sets
\begin{equation*}
    \bad_{A}(\varepsilon):=\left\{\bx \in \R^{m}: \left\| A\bq -\bx -\bp\right\|_{\infty} \geq \varepsilon \|\bq\|_{\infty}^{-\tfrac{n}{m}} \quad \text{ for all }\, (\bp,\bq) \in \Z^{m}\times (\Z^{n}\setminus \{\mathbf{0}\}) \right\},
\end{equation*}
and \begin{equation*}
    \bad_{A}:=\bigcup_{\varepsilon>0} \bad_{A}(\varepsilon).
\end{equation*}

The set $\bad_{A}$, often called the set of twisted badly approximable points, is a classical object of study. The term ``twisted'' comes from the fact we are now fixing $A\in \R^{m\times n}$ and considering points $\bx$ that avoid a shrinking neighbourhood of the orbit of $(A\bq \mod \Z^{m})_{\bq\in \Z^{n}}$, rather than the classical notion of badly approximable points where the $\bx$ is fixed ($\bx=\mathbf{0}$ in the homogeneous set of badly approximable points) and $A$ varies.

We recall the following \textit{non-singular} condition. An $m\times n$ matrix $A$ is \textit{singular} if for any $\varepsilon>0$ there exists $H_{0}$ such that for all $H>H_{0}$ there exists solution $(\bp,\bq)\in \Z^{m}\times (\Z^{n}\setminus \{\mathbf{0}\})$ to the system
\begin{equation*}
    \| A\bq-\bp\|_{\infty}<\varepsilon H^{-\frac{n}{m}}, \quad 1\leq \|\bq \|_{\infty} \leq H.
\end{equation*}
Otherwise a matrix $A\in \R^{m\times n}$ is said to be non-singular. Moreover, a matrix $A\in \R^{m\times n}$ is said to be \textit{singular-on-average} if
\begin{equation*}
    \lim_{M\to\infty} \frac{1}{M}\# \left\{ 1\leq \ell \leq M: \begin{array}{c} \text{there is a solution $(\bp,\bq)\in \Z^{m}\times (\Z^{n}\setminus \{\mathbf{0}\})$ to}\\
     \|A\bq-\bp \|_{\infty}<\delta 2^{-\ell \tfrac{n}{m}}, \quad 1\leq \|\bq\|_{\infty}< 2^{\ell} 
     \end{array}\right\}=1,
\end{equation*}
for every $\delta>0$. Clearly a matrix that is not singular-on-average, so there exists $\delta>0$ such that the above lower limit is strictly less than one, is non-singular.

For $m=n=1$ it was shown in \cite{DHKim07} that for any irrational $\alpha\in \R$ the set $\bad_{\alpha}$ is a Lebesgue nullset. It was shown in \cite[Theorem 1.6]{MoshRaoShap24} that for any $A\in \R^{m\times n}$ non-singular the set $\bad_{A}$ is a Lebesgue nullset, see also \cite{BDGWTori}. This left open the question of the Lebesgue measure of $\bad_{A}$ for general $A$ which has been recently resolved, see \cite[Proposition 2.2]{KimNec26}. They show that for any $\varepsilon>0$ we have $\lambda_{m}(\bad_{A}(\varepsilon))=\lambda_{m}(\bad_{A})=0$ if and only if
\begin{equation*}
    \sum_{k=1}^{\infty} \left(\#\left\{(\bp,\bq)\in \Z^{m}\times \Z^{n}: \|A\bq-\bp\|_{\infty}\leq r=2^{-k\frac{n}{m}}\, \|\bq\|_{\infty}\leq 2^{k} \right\} \right)^{-1}=\infty.
\end{equation*}
In particular, \cite[Corollary 2.4]{KimNec26} recovers \cite[Theorem 1.6]{MoshRaoShap24}. This naturally leads one to studying the dimension of $\bad_{A}$ and $\bad_{A}(\varepsilon)$ for $A$ satisfying the above divergence condition.

It was shown in \cite{BHKV10} that for any $A\in \R^{m\times n}$
\begin{equation*}
    \dimh \bad_{A}=m.
\end{equation*}
See \cite{Tseng09, EinTse11, BenMos17} for various generalisations of the above statements. 

The dimension of the smaller sets $\bad_{A}(\varepsilon)$ has also been studied. Lim, de Saxce and Shapira \cite[Theorem 1.5]{LdSS19} proved that for any $\varepsilon>0$ there exists $\delta>0$ such that for Lebesgue almost all $A\in  \R^{m}$, 
\begin{equation*}
\dimh \bad_{A}(\varepsilon)<m-\delta.    
\end{equation*}
Conversely, in \cite{BugKLR21} it was shown that for $m=n=1$ the set $\bad_{\alpha}(\varepsilon)$ has full Hausdorff dimension for some $\varepsilon>0$ if and only if $\alpha$ is singular-on-average. This result was later generalised to higher dimensions in \cite[Theorem 1.3]{KimKimLim25}.

In \cite[Theorem 1.2]{KimKimLim25} it was shown that for any $A\in \R^{m\times n}$ not singular-on-average there exists a constant $c(A)>0$ such that for any $\varepsilon>0$,
\begin{equation*}
    \dimh \bad_{A}(\varepsilon)\leq m - c(A)\frac{\varepsilon^{m}}{|\log \varepsilon|},
\end{equation*}
thus improving upon the result of \cite{LdSS19}. We prove the following result

\begin{theorem} \label{bad irrational}
Let $A\in \R^{m\times n}$ be non singular, $\delta>0$ and define
\begin{equation*}
    L_{A}(\delta):=\left\{ \ell \in \N : \begin{array}{c} \text{there is \textbf{no} solution $(\bp,\bq)\in \Z^{m}\times (\Z^{n}\setminus \{\mathbf{0}\})$ to}\\
     \|A\bq-\bp \|_{\infty}<\delta 2^{-\ell \tfrac{n}{m}}, \quad 1\leq \|\bq\|_{\infty}< 2^{\ell} 
     \end{array}\right\}.
\end{equation*}
Then, for $\varepsilon>0$ sufficiently small depending on $\delta>0$, we have that
\begin{equation*}
    \dimp \bad_{A}(\varepsilon) \leq m - C_{A} \varepsilon^{m},
\end{equation*}
with constant 
\begin{equation*}
    C_{A}=C\left(\frac{2}{\delta}\right)^{m}(\delta^{-m}+1)^{-(m+n)} \liminf_{M\to\infty} \frac{\# (L_{A}(\delta)\cap [1,M])}{M},
\end{equation*}
for $C>0$ depending on $n$ and $m$ only.

\end{theorem}

\begin{remark} \rm
Due to the aforementioned results of \cite{BugKLR21,KimKimLim25, KimNec26} the dependence on the Diophantine properties of $A$ is crucial. The techniques of the previously mentioned results are significantly different from that presented here. In particular, many of the papers translate it into a problem in homogeneous dynamics using the \textit{Dani Correspondence} \cite{Dani85}. 
\end{remark}

\begin{remark}
    Note that the dimension bound may not necessarily be non-trivial (full). In particular, it may be possible for a matrix $A$ to be non singular which implies there exists some $\delta>0$ for which
    \begin{equation*}
        \limsup_{M\to\infty}\frac{\#(L_{A}(\delta)\cap [1,M])}{M}>0,
    \end{equation*}
    yet
    \begin{equation*}
    \liminf_{M\to\infty}\frac{\#(L_{A}(\delta)\cap [1,M])}{M}=0.
    \end{equation*}
\end{remark}

As a corollary to Theorem~\ref{bad irrational} one can obtain the following result.

\begin{corollary}
Let $A\in \R^{m\times n}$ be a badly approximable matrix with best approximation constant $c(A)>0$. That is, 
\begin{equation} \label{eq: bad def}
    \|A\bq-\bp\|_{\infty}\geq c(A) \|\bq\|_{\infty}^{-\frac{n}{m}} \quad \text{ for all pairs } \, (\bp,\bq)\in \Z^{m}\times (\Z^{n}\setminus \{\mathbf{0}\}).
\end{equation}
Then, for $\varepsilon>0$ sufficiently small we have that
\begin{equation*}
    \dimp \bad_{A}(\varepsilon) \leq m-C_{A} \varepsilon^{m},
\end{equation*}
with $C_{A}=C\left(\tfrac{2}{c(A)}\right)^{m}(c(A)^{-m}+1)^{-(m+n)}$.
\end{corollary}
\begin{proof}
Suppose $\ell\in \N$ and  $\ell \not\in L_{A}(c(A))$. Then there exists solution $(\bp,\bq)\in \Z^{m}\times (\Z^{n}\setminus \{\mathbf{0}\})$ solving
\begin{equation*}
    \|A\bq-\bp\|_{\infty}<c(A) 2^{-\ell \frac{n}{m}}, \quad 1\leq \|\bq\|_{\infty}<2^{\ell}.
\end{equation*}
This implies 
\begin{equation*}
    \|A\bq-\bp\|_{\infty}<c(A)\|\bq\|_{\infty}^{-\frac{n}{m}},
\end{equation*}
which contradicts \eqref{eq: bad def}. Hence $\ell \in L_{A}(c(A))$ and so $L_{A}(c(A))=\N$. Looking at the lower limit of the frequency of $L_{A}(c(A))$ and using Theorem~\ref{bad irrational} we obtain the result.
\end{proof}

\subsubsection{Proof of Theorem~\ref{bad irrational}}
 The following transference principle between homogeneous and inhomogeneous Diophantine approximation, see \cite[Chpt. V Theorem VI]{Cassels}, can be used to obtain the following crucial proposition:

\begin{proposition}[{\cite[Corollary 3.4]{TKim24}}] \label{prop: Kim result}
    Let $A\in \R^{m\times n}$ be non-singular. Then for any $\bx \in \R^{m}$ and any $\ell \in L_{A}(\delta)$, there exists solution $(\bp,\bq) \in \Z^{m}\times \Z^{n}$ to the system of inequalities
    \begin{equation*}
    \|A\bq -\bx-\bp \|_{\infty}\leq \tfrac{1}{2} (\delta^{-m}+1)\delta 2^{-\ell \tfrac{n}{m}}, \quad 0\leq \| \bq \|_{\infty}\leq \tfrac{1}{2}(\delta^{-m}+1) 2^{\ell}.
    \end{equation*}
\end{proposition}

Armed with this result we can prove Theorem~\ref{bad irrational}

\begin{proof}[Proof of Theorem~\ref{bad irrational}]
Firstly, note that if
\begin{equation*}
    F(L_{A}(\delta)):= \liminf_{M\to\infty} \frac{\# (L_{A}(\delta)\cap [1,M])}{M}=0,
\end{equation*}
then there is nothing to prove, so we may assume that $F(L_{A}(\delta))>0$, hence $L_{A}(\delta)$ is infinite.
In the notation of \S~\ref{sec: main result} let 
\begin{align*}
    J=\Z^{n}\setminus \{\mathbf{0}\}, \quad R_{\bq}=\{A\bq +\bp: \bp \in \Z^{m}\}, \quad
    h(\bq)&=\|\bq\|_{\infty},
\end{align*}
so $\cR=\left\{A\bq +\bp: \bq \in \Z^{n}\setminus \{\mathbf{0}\} \right\}\subset \R^{m}$. Note that
\begin{equation*}
    \bad_{A}(\varepsilon)\subset \left(\bigcup_{1\leq i \leq m}\{(x_{1},\ldots, x_{m})\in \R^{m}: x_{i}\in \Z\}\right) \cup \left(\bigcup_{\bz \in \Z^{m}}\cB_{\cR}\left(\varepsilon, \frac{n}{m}, (0,1)^{m}+\bz\right) \right)
\end{equation*}
by definition. The first part of the above cover is all integer translates of coordinate hyperplanes, which has packing dimension equal to $m-1$. 
To apply Theorem~\ref{main} we need $\|\bq\|_{\infty}\geq 1$, which Proposition~\ref{prop: Kim result} does not currently supply. Fix $\bz \in \Z^{m}$. Let $(B_{i})_{i\in\N}$ be a collection of open balls such that $\overline{B_{i}}\subset (0,1)^{m}+\bz$ and $$\bigcup_{i\in \N} B_{i}=(0,1)^{m}+\bz.$$
Fix $i \in \N$. For any $\bx \in B_{i}\subset (0,1)^{m}+\bz$ we have that $\|\bx-\bp\|_{\infty}>\dist(B_{i}, \left((0,1)^{m}+\bz\right)^{c})=e_{i}>0$ for all $\bp\in \Z^{m}$. Therefore, for $\ell\geq \ell_{0,i}$ sufficiently large depending on $\delta$ and $i$, by Proposition~\ref{prop: Kim result} we have that for any $\bx \in B_{i}$ and any $\ell \in L_{A}(\delta)_{\geq \ell_{0,i}}$ there exists solution $(\bp,\bq)\in \Z^{m}\times \Z^{n}$ to the system
\begin{equation*}
    \|A\bq-\bx-\bp\|_{\infty}<\frac{1}{2}(\delta^{-m}+1)\delta 2^{-\ell\frac{n}{m}}, \quad 1\leq \|\bq\|_{\infty}<(\delta^{-m}+1)2^{\ell}.
\end{equation*}
So condition \eqref{eq:dirichlet type cond} is satisfied for $\cU=B_{i}$ with 
\begin{equation*}
    a=0, \qquad b=\frac{n}{m}, \qquad c_{1}=\frac{1}{2}(\delta^{-m}+1)\delta, \qquad c_{2}=(\delta^{-m}+1),
\end{equation*}
and $\cS=L_{A}(\delta)_{\geq \ell_{0,i}}$. Applying Theorem~\ref{main} to $\cB_{\cR}\left(\varepsilon, \frac{n}{m},B_{i}\right)$ we see that
\begin{equation*}
    \dimp \cB_{\cR}\left(\varepsilon, \frac{n}{m},B_{i}\right)\leq  m- C_{A,i} \varepsilon^{m},
\end{equation*}
for 
\begin{align*}
    C_{A,i}&=C\left(\frac{2}{\delta}\right)^{m}(\delta^{-m}+1)^{-(m+n)}\liminf_{M\to \infty} \frac{\#(L_{A}(\delta)_{\geq\ell_{0,i}}\cap [1,M])}{M} \\
    &=C\left(\frac{2}{\delta}\right)^{m}(\delta^{-m}+1)^{-(m+n)} F(L_{A}(\delta)). 
\end{align*}
Clearly $C_{A,i}$ is independent of $i$, so let $C_{A}=C_{A,i}$. By countable stability of the packing dimension, first over $i \in \N$ then over $\bz \in \Z^{m}$, we have that
\begin{equation*}
    \dimp \bad_{A}(\varepsilon) \leq \max \left\{ m-1, m-C_{A}\varepsilon^{m} \right\}= m-C_{A}\varepsilon^{m}
\end{equation*}
for $\varepsilon>0$ sufficiently small.
\end{proof}

\bibliographystyle{abbrv}
\bibliography{biblio2}

\end{document}